\documentclass[draft, 12pt]{amsart}
\usepackage[english]{babel}

\usepackage{amssymb} 
\usepackage{amsmath}
\usepackage[all,cmtip]{xy}
\usepackage[pagewise]{lineno}

\usepackage{color}
\usepackage{graphicx}
\usepackage{tikz}
\usepackage{circuitikz}
\usetikzlibrary{arrows.meta}

\newtheorem{thm}{Theorem}
\newtheorem{defi}{Definition}
\newtheorem{prop}[thm]{Proposition}
\newtheorem{cor}[thm]{Corollary}
\newtheorem{lem}[thm]{Lemma}

\usepackage{ifdraft, blindtext, verbatim}
\ifdraft{%
    \usepackage{todonotes}}{ 
    \usepackage[disable]{todonotes} 
}

\begin{document}
\title[Weinstein's relative morphism]
{A relative version of Weinstein's morphism
}

\author{Andr\'es Pedroza}  
\address{Facultad de Ciencias\\
           Universidad de Colima\\
           Bernal D\'{\i}az del Castillo No. 340\\
           Colima, Col., Mexico 28045}
\email{andres\_pedroza@ucol.mx}

\author{Daniel Villanueva-Vega}  
\address{Department of Mathematics\\
           Notre Dame University}
\email{dvillanu@nd.edu}

 \begin{abstract}
 We introduce a relative version of Weinstein's morphism associated to a
 Lagrangian submanifold \(L\) of a closed symplectic manifold
 \((M,\omega)\). We construct a homomorphism
 \[
 \mathcal{A}^{(L,\alpha)}\colon
 \pi_{1}\bigl(\operatorname{Ham}_0(M,L)\bigr)
 \longrightarrow
 \mathbb{R}/\mathcal{P}_{2}(M,\omega),
 \]
 and establish its basic properties. We compute the invariant
 in several examples, including latitude circles in \(S^2\),
and more generally, Lagrangian fibers
in toric manifolds. 
These computations provide
 obstructions to Ha\-mil\-to\-nian isotopy of Lagrangian submanifolds and detect
 nontrivial, and in some cases infinite-order, elements in the homotopy
 groups of their Hamiltonian stabilizers.
 \end{abstract}

\keywords{Weinstein's morphism, Lagrangian submanifolds, Hamiltonian diffeomorphism}
\subjclass{Primary: 57S05, 53D35 Secondary:  57R17, 53D05}
\maketitle

\section{Introduction}


The topology of the group of Hamiltonian diffeomorphisms of a symplectic manifold is difficult to determine. The presence of a Lagrangian submanifold leads to a related, but distinct, question: what can be said about Hamiltonian isotopies that preserve the Lagrangian at every time?  Let $(M,\omega)$ be a closed symplectic manifold and let $L\subset M$ be a closed Lagrangian submanifold. We write
\[\text{Ham}(M,L) = \{\phi\in \text{Ham}(M,\omega):\phi(L)=L\}
\]
and denote its identity component by $\text{Ham}_0(M,L)$. The inclusion \[\text{Ham}_0(M,L)\hookrightarrow\text{Ham}(M,\omega)\]
induces a homomorphism on fundamental groups.  Describing these relative homotopy classes, and finding obstructions are natural problems in the topology of the group of Hamiltonian diffeomorphisms.

One classical approach to the topology of symplectic transformation groups is through homomorphisms defined by symplectic action or flux. Weinstein's action morphism associates to a Hamiltonian loop an element of the quotient of $\mathbb{R}$ by the spherical period group \cite{weinstein-coho}. Another approach is provided by the relative Seidel morphism of S. Hu and F. Lalonde, defined using Lagrangian Floer theory, and its applications to homological Lagrangian monodromy \cite{hulalonde}. 
Also,  A. S. Demir  defined in \cite{Demir} 
  a relative version of  Calabi morphism on the universal cover of  ${\textup{Ham}(M,L)}$.
More recently, J. Brendel, J.-P. Chasse, and R. Leclercq introduce a Lagrangian flux-monodromy morphism on the fundamental group of a space of Lagrangian submanifolds \cite{BCL}.  Also, Y. Ozan studies in  \cite {Ozan} the larger group
 $\textup{Symp}(M,L)$ by extending the definition of the flux morphism.

In this paper we consider an elementary relative version of Weinstein's morphism, based on the symplectic areas swept out by points of a Lagrangian under Hamiltonian loops preserving it. Our aim is to make this area construction explicit, establish its naturality, and use it to obtain obstructions without invoking Floer homology, quantum homology, or pseudoholomorphic-curve theory. Throughout the construction we assume that $M$ is simply connected, and we denote by \[\mathcal{P}_2(M,\omega) = \Big\{ \int_A \omega: A \in H_2 (M;\mathbb{Z})\Big\}\]
the spherical period group. For an oriented Lagrangian $L$ endowed with a volume form $\alpha$, we define the morphism \[\mathcal{A}^{(L,\alpha)}: \pi_1\left(\text{Ham}_{0}(M,L)\right) \rightarrow \mathbb{R}/\mathcal{P}_2(M,\omega).\] When $L$ is connected, it has a simpler description. If $\{\psi_t\}_{t\in[0,1]}$ is a loop in $\text{Ham}_0(M,L)$ and $p\in L$, the evaluation loop $t\mapsto\psi_t(p)$ bounds a disk $u_p: D^2 \rightarrow M,$ then \[ \mathcal{A}^{(L,\alpha)}([\psi_t]) = \int_{D^2}u_{p}^{*}\omega 
\quad\text{in }\mathbb{R}/\mathcal{P}_2(M,\omega).
\]
Because the symplectic form vanishes on $L$, the right-hand side is independent of the capping disk and the point $p$, modulo spherical periods. In particular, for connected Lagrangians the homomorphism does not depend on the auxiliary volume form. Equivalently, it is the composition of the evaluation in $L$ with the relative symplectic-area character. This elementary factorization is useful for both calculations and understanding the limits of the construction. 

The first results establish the homomorphism property and naturality under Hamiltonian conjugation. Consequently, the image of the morphism is an invariant of the Hamilotnian isotopy class of the Lagrangian. Comparing these images gives a sufficient obstruction to Hamiltonian isotopy. 

In particular, we use the relative Weinstein morphism to detect nontrivial elements, including elements of infinite order, in the fundamental group $\pi_{1}\left( \text{Ham}_{0}(M,L)\right)$, and to recover the fact that certain Lagrangian submanifolds are not Hamiltonian isotopic.

For a connected Lagrangian, the same description also yields a fixed-point obstruction: if $\mathcal{A}^{(L,\alpha)}([\psi_t])\neq 0$, no representative of that relative homotopy class can fix a point of $L$ throughout the loop. These conclusions require only the Lagrangian condition and the elementary properties of the symplectic area.

\section{Main results}

We begin with the simplest example, the $2$-sphere. Although elementary, this example already exhibits an important observation: a Hamiltonian loop that has finite order in the absolute Hamiltonian group can have infinite order in the relative group.

\begin{thm}
\label{t:sphere}
Let \((S^2,\omega)\) be the unit sphere, and for \(h\in[0,1)\) let
\[
L_h
=
\{(x,y,h)\in S^2\mid x^2+y^2+h^2=1\}.
\]
Let \(\{\psi_t\}\) be the Hamiltonian loop given by rotation about the
\(z\)-axis. If
${(1-h)}/{2}$ is irrational
then
\[
[\psi_t]\in
\pi_1\bigl(\operatorname{Ham}_0(S^2,L_h)\bigr)
\]
has infinite order.
\end{thm}

Recall that $\operatorname{Ham}(S^2,\omega)$ deformation retracts onto
$SO(3)$, thus
$\pi_1\bigl(\operatorname{Ham}(S^2,\omega)\bigr)\cong\mathbb{Z}_2$,
and $[\psi_t]$ is its generator. \\

A family of higher-dimensional examples comes from symplectic toric manifolds. For interior points in the Delzant polytope $\Delta(M,\omega)$ of a toric manifold $(M,\omega)$, let $L_{\bf a} = \mu^{-1}({\bf a})$ be an interior Lagrangian torus fiber. Every integral vector $\xi \in \mathbb{Z}^n$ gives a circle action 
$\{\psi^{\xi}_t\}$
preserving $L_a$. If $v$ is a vertex of the Delzant polytope, an elementary capping construction gives 
\begin{equation*}\label{eqn:affine}
    \mathcal{A}^{(L,\alpha)}([\psi^{\xi}_t]) =  \langle a -v, \xi \rangle  
    \quad\text{in }\mathbb{R}/\mathcal{P}_2(M,\omega).
\end{equation*}
For a primitive inward face normal, this is the integral-affine distance from $a$ to the corresponding facet, modulo spherical periods.

We illustrate this in the case of
\((\mathbb{C}P^n,\omega_{\mathrm{FS}})\),
 with the symplectic form normalized
so that the symplectic area of a projective line is \(1\), which is the
standard normalization in the toric setting. Henceforth,
$\mathcal{P}_{2}  (\mathbb{C}P^n,\omega_{\mathrm{FS}})=\mathbb{Z}$.
Let \(\Delta (\mathbb{C}P^n,\omega_{\mathrm{FS}}) \subset  \mathfrak{t}^*  \) denote the
corresponding moment polytope, and for \(1\leq k\leq n\), let
\(\{\psi_t^{(k)}\}\) denote the Hamiltonian circle action corresponding to
the \(k\)-th component of the moment map $\mu\colon (\mathbb{C}P^n,\omega_{\mathrm{FS}})
\to \mathfrak{t}^*$,
\begin{eqnarray*}
\psi^{(k)}_t
[z_0,\ldots,z_n]
=
[z_0,\ldots,e^{2\pi i t}z_k,\ldots,z_n].
\end{eqnarray*}

For
$
\mathbf a=(a_1,\ldots,a_n)\in\operatorname{int}(\Delta (\mathbb{C}P^n,\omega_{\mathrm{FS}}) ),
$
 let \(L_{\mathbf a}=\mu^{-1}(\mathbf a)\) be the corresponding
Lagrangian torus.  Then $\{ \psi^{(k)}_t \} $
 is a loop in $\textup{Ham}_0(\mathbb{C}P^n, L_{\mathbf a}).$

\begin{thm}
\label{t:aoncpn}
For every \(1\leq k\leq n\) and any volume form $\alpha$ in $L_{\mathbf a}$,
\[
\mathcal{A}^{(L_{\mathbf a},\alpha)}
([\psi^{(k)}_t])
=
a_k
\quad\text{in }\mathbb{R}/\mathbb{Z}.
\]
In particular if $a_k$ is irrational, then 
\begin{eqnarray*}
[\psi^{(k)}_t]\in \pi_1(\textup{Ham}_0(\mathbb{C}P^n, L_{\mathbf a}))
\end{eqnarray*}
has infinite order.
\end{thm}

Observe that the second part of this result does not apply to the Clifford
torus $L_{(1/(n+1),\ldots,1/(n+1))}\subset\mathbb{C}P^n$; for $n=1$, this is
the equator in $S^2$, to which Theorem~\ref{t:sphere} does not apply either.

A consequence of Theorem~\ref{t:aoncpn} is that it can be used to
distinguish Lagrangian fibers in
\((\mathbb{C}P^n,\omega_{\mathrm{FS}})\) that are not Hamiltonian isotopic.
Although this is a well-known fact, our approach provides a direct
application of the relative Weinstein morphism.

\begin{cor}
\label{c:toric-image}
For  $\mathbf a=(a_1,\ldots,a_n)\in\operatorname{int}(\Delta (\mathbb{C}P^n,\omega_{\mathrm{FS}})   )$ and  any volume form $\alpha$ on
$L_{\mathbf a}$,
\[
\operatorname{Im}\mathcal{A}^{(L_{\mathbf a},\alpha)}
= \mathbb{Z}\langle 1, a_1,\ldots ,a_n \rangle  /\mathbb{Z}
\subset\mathbb{R}/\mathbb{Z}.
\]
Consequently, if $\mathbf b=(b_1,\ldots,b_n)\in\operatorname{int}(\Delta(\mathbb{C}P^n,\omega_{\mathrm{FS}})   )$
satisfies
\[
 \mathbb{Z}\langle 1, a_1,\ldots ,a_n \rangle 
\neq
 \mathbb{Z}\langle 1, b_1,\ldots ,b_n \rangle 
\]
as subgroups of $\mathbb{R}$, then $L_{\mathbf a}$ and $L_{\mathbf b}$ are
not Hamiltonian isotopic in $(\mathbb{C}P^n,\omega_{\mathrm{FS}}).$
\end{cor}

In Section~\ref{sec:toric}, we derive an  analogous result for the case of toric manifolds.

As mentioned above, in the case when the Lagrangian is connected the
morphism $\mathcal{A}^{(L,\alpha)}$ admits a simple description.
This alternative description of the  morphism  provides information
about fixed points of Hamiltonian loops on the Lagrangian submanifold.
If
$\mathcal{A}^{(L,\alpha)}([\psi_t])\neq0$, then no loop homotopic to
$\{\psi_t\}$ in $\operatorname{Ham}_0(M,L)$ has a point of $L$ that is fixed
for all $t$; in fact, the orbit of every point of $L$ under such a loop is
non-contractible in $L$. Since $M$ is simply
connected, no such restriction holds for orbits of points outside $L$.



\begin{thm}
\label{t:nonfix}
Let $(M,\omega)$ be a closed, simply connected symplectic manifold, and let
$L\subset M$ be a connected, oriented Lagrangian submanifold with volume
form $\alpha$. If $\{\psi_t\}$ is a loop in $\operatorname{Ham}_0(M,L)$
based at the identity such that
$\mathcal{A}^{(L,\alpha)}[\psi_t]\neq0$, then for every loop
$\{\psi'_t\}$ homotopic to $\{\psi_t\}$ in $\operatorname{Ham}_0(M,L)$
there is no point $p\in L$ with $\psi'_t(p)=p$ for all $t\in[0,1]$.
\end{thm}




\section{Definition of the morphism $\mathcal{A}^{(L, \alpha)} $}

In this section, we define the morphism
\[
\mathcal{A}^{(L, \alpha)} \colon \pi_{1}(\textup{Ham}_0(M, L)) \to \mathbb{R}/
\mathcal{P}_{2}(M, \omega),
\]
which requires showing that the resulting construction is independent of the
various auxiliary choices made in its definition.
Throughout this section, \( (M,\omega) \) denotes a simply connected, closed symplectic manifold of 
dimension \(2n\), and \(L\subset (M,\omega)\)  a Lagrangian submanifold.
As will be clear from the definition that  the assumption that \( (M, \omega) \) is closed can be 
relaxed. In the non-compact case, however, one must consider compactly supported Hamiltonian diffeomorphisms.

Recall that the period group
\(\mathcal{P}_{2}(M,\omega)\subset\mathbb{R}\) is defined as the image of
the pairing
\[
\langle[\omega],\cdot\rangle\colon
H_{2}(M;\mathbb{Z})\longrightarrow\mathbb{R}.
\]
This group was used  by Weinstein in \cite{genwen} in the definition of
the morphism
$
\mathcal{A}\colon
\pi_1(\operatorname{Ham}(M,\omega))
\longrightarrow
\mathbb{R}/\mathcal{P}_2(M,\omega).
$
Note that since $M$ is simply connected, the Hurewicz map
$\pi_2(M)\to H_2(M;\mathbb{Z})$ is an isomorphism; in particular, every
class in $H_2(M;\mathbb{Z})$ is spherical.

Define the map
$
A\colon C^\infty(S^{1},M)\longrightarrow
\mathbb{R}/\mathcal{P}_{2}(M,\omega)
$
as follows:
\[
A(\gamma):=
\int_{D^{2}}u_\gamma^*(\omega),
\]
where \(\gamma\in C^\infty(S^{1},M)\) and
\(u_\gamma\colon D^{2}\to M\) is a smooth map whose restriction to
\(\partial D^{2}\) agrees with \(\gamma\). 
Since \(A\) takes values in
\(\mathbb{R}/\mathcal{P}_{2}(M,\omega)\), its definition is independent of the
choice of capping \(u_\gamma\).

For \(p\in L\) and
\(\{\psi_t\}\in
\Omega^{1}(\operatorname{Ham}_0(M,L))\),
define $
\gamma(\{\psi_t\},p)\colon
S^{1}\longrightarrow L$  by  $ \gamma(\{\psi_t\},p)(t)
=\psi_t(p).$
We denote by
$
C(\{\psi_{t}\},p)\subset L
$
the trace of \(p\) under the family
\(\{\psi_t\}_{{t}\in S^{1}}\), that is, the image of
\(\gamma(\{\psi_{t}\},p)\).

\begin{lem}
\label{l:lemhomo}
Fix \(p\in L\). If
$
[\psi^{(0)}_{t}]
=
[\psi^{(1)}_{t}]
\in \pi_{1}(\operatorname{Ham}_0(M,L)),
$
then
\[
A\bigl( \gamma(\{\psi^{(0)}_{t}\},p)  \bigr)
=
A\bigl(  
\gamma(\{\psi^{(1)}_{t}\},p) \bigr)
\]
in \(\mathbb{R}/\mathcal{P}_{2}(M,\omega)\).
\end{lem}
\begin{proof}
Consider a homotopy 
\[
\boldsymbol{\Psi} \colon [0,1] \times S^{1} \to \textup{Ham}_0(M, L)
\]
between  \(\{\psi^{(0)}_{t}\}\) and \(\{\psi^{(1)}_{t}\}\).  
Evaluating \(\boldsymbol{\Psi}\) at the base point \(p \in L\) induces a homotopy in \(M\),
\[
\boldsymbol{\Psi}_p \colon [0,1] \times S^{1} \to M,
\]
between the traced circles \(C(\{\psi^{(0)}_{t}\},p)\) and
\(C(\{\psi^{(1)}_{t}\},p)\).
Moreover, since each \(\psi^{(s)}_{t}\) preserves \(L\), the image of the  homotopy \(\boldsymbol{\Psi}_p\) 
lies inside the Lagrangian submanifold \(L \subset (M, \omega)\).

Thus, the caps 
$ u_{  \gamma(\{\psi^{(0)}_{t}\} ,p) }, u_{  \gamma(\{\psi^{(1)}_{t}\} ,p) } 
\colon (D^2,\partial D^2)\longrightarrow(M,L)$
together with the image of the  homotopy  \(\boldsymbol{\Psi}_p([0,1] \times S^{1}) \subset L\), form a  \(2\)-sphere in \(M\).  
Therefore, by definition of the group \(\mathcal{P}_{2}(M, \omega)\) and the fact that $\omega$ vanishes  on $L$,
 we have
\[
\int_{D^{2}} u_{  \gamma(\{\psi^{(0)}_{t}\} ,p) }^*(\omega) =
 \int_{D^{2}} u_{ \gamma(\{\psi^{(1)}_{t}\},p) }^*(\omega)
\]
in \(\mathbb{R}/\mathcal{P}_{2}(M, \omega)\).
\end{proof}

For \(\{\psi_{t}\}\in
\Omega^{1}(\operatorname{Ham}_0(M,L))\),   define the map
$\operatorname{ev}_{\psi_{t}}\colon M\longrightarrow C^\infty(S^{1},M)$
by $ \operatorname{ev}_{\psi_{t}}(p)(t) :=
\psi_{t}(p).$ Composing with \(A\), we obtain the map
\[
A\circ\operatorname{ev}_{\psi_{t}}\colon
M\longrightarrow
\mathbb{R}/\mathcal{P}_{2}(M,\omega).
\]
In order to integrate this map, we need to lift it to a
real-valued function. To this end, note that
$ \operatorname{ev}_{\psi_{t}}^*(dA)$
is a closed
\(\mathbb{R}/\mathcal{P}_{2}(M,\omega)\)-valued \(1\)-form on \(M\).
Since \(M\) is simply connected, there exists a function
$ F_{\psi_{t}}\colon M\longrightarrow\mathbb{R}$
such that
\[
dF_{\psi_{t}}
=
\operatorname{ev}_{\psi_{t}}^*(dA).
\]
By adding a suitable constant to \(F_{\psi_{t}}\), we may assume
that \(F_{\psi_{t}}\) is a lift of
\(A\circ\operatorname{ev}_{\psi_{t}}\). Thus,
\[
F_{\psi_{t}}(p)
-
A\bigl(\psi_{t}(p)\bigr)
=
0
\quad\text{in }\mathbb{R}/\mathcal{P}_{2}(M,\omega).
\]

Finally, 
 we require the
Lagrangian submanifold \(L\) to be oriented and endowed with a fixed volume form
\(\alpha\).
Note that, depending on the context, we will regard \(L\) either as a
Lagrangian submanifold of \((M,\omega)\) or as an abstract manifold equipped
with a Lagrangian embedding
$
\iota\colon L\longrightarrow (M,\omega).
$
The intended meaning will be clear from the context.

Given
$
\{\psi_{t}\}\in
\Omega^{1}(\operatorname{Ham}_0(M,L)),
$
define
\begin{eqnarray*}
\mathcal{A}^{(L,\alpha)}
\bigl(\{\psi_{t}\}\bigr)
&:=&
\frac{1}{\operatorname{Vol}(L,\alpha)}
\int_L F_{\psi_{t}}\,\alpha \\
&=&
\frac{1}{\operatorname{Vol}(L,\alpha)}
\int_L
\left(
\int_{D^{2}}
 u_{ \gamma(\{\psi_{t}\},p) }^*(\omega)
\right)\alpha(p)
\end{eqnarray*}
as an element of
$
\mathbb{R}/\mathcal{P}_{2}(M,\omega).
$
By Lemma \ref{l:lemhomo}, the value of
\(\mathcal{A}^{(L,\alpha)}\) is invariant under homotopies of the loop
\(\{\psi_{t}\}\).

\begin{defi}
Let \((M,\omega)\) be a simply connected, compact symplectic manifold of
dimension \(2n\).
Let \(L\) be an oriented \(n\)-dimensional manifold equipped with a fixed
volume form \(\alpha\), and suppose that \(L\) is a Lagrangian submanifold of
\((M,\omega)\). We define
\[
\mathcal{A}^{(L,\alpha)}\colon
\pi_{1}(\operatorname{Ham}_0(M,L))
\longrightarrow
\mathbb{R}/\mathcal{P}_{2}(M,\omega)
\]
by
\[
\mathcal{A}^{(L,\alpha)}
\bigl([\psi_{t}]\bigr)
:=
\frac{1}{\operatorname{Vol}(L,\alpha)}
\int_L
\left(
\int_{D^{2}}
u_{ \gamma(\{\psi_{t}\},p) }^*(\omega)
\right)\alpha(p).
\]
\end{defi}

Geometrically, the map \(\mathcal{A}^{(L,\alpha)}\) can be interpreted as
the average, over the Lagrangian submanifold \(L\), of the symplectic
\(2\)-volume of a disk bounding the trace of a point \(p\in L\) under the
family \(\{\psi_{t}\}\), where the average is understood modulo
\(\mathcal{P}_{2}(M,\omega)\).

Note that if the Lagrangian submanifold \(L\subset(M,\omega)\) satisfies
$\pi_{1}(L)=0,$
then
$\mathcal{A}^{(L,\alpha)}\equiv 0.$
Also, if
$
[\psi_{t}]\in
\pi_{1}(\operatorname{Ham}_0(M,L))
$
is represented by a family
\(\{\psi_{t}\}_{t\in S^{1}}\)
 such that each
\(\psi_{t}\) fixes \(L\) pointwise, then
$
\mathcal{A}^{(L,\alpha)}
\bigl([\psi_{t}]\bigr)=0.$
In the sections that follow, we compute
\(\mathcal{A}^{(L,\alpha)}\) for several classes of symplectic
manifolds and Lagrangian submanifolds. These examples illustrate the
geometric and topological information captured by the relative Weinstein
morphism.

The discussion so far has focused on a Lagrangian submanifold
\(L\subset(M,\omega)\), motivated by the fundamental role that Lagrangian
submanifolds play in symplectic topology. However, the same construction
applies more generally to an oriented isotropic submanifold
\(L'\subset(M,\omega)\). Indeed, following the same argument as above, one
obtains a map
\[
\pi_{1}(\operatorname{Ham}_0(M,L'))
\longrightarrow
\mathbb{R}/\mathcal{P}_{2}(M,\omega).
\]

An important feature of the map
\(\mathcal{A}^{(L,\alpha)}\) is that it is a group homomorphism.
This structure allows the morphism to detect nontrivial topological
information about the group \(\operatorname{Ham}_0(M,L)\).

\begin{prop}
The map
\[
\mathcal{A}^{(L,\alpha)}\colon
\pi_{1}(\operatorname{Ham}_0(M,L))
\longrightarrow
\mathbb{R}/\mathcal{P}_{2}(M,\omega)
\]
is a homomorphism of abelian groups.
\end{prop}

\begin{proof}
Let
$
[\psi^{(0)}_{t}],
[\psi^{(1)}_{t}]
\in
\pi_{1}(\operatorname{Ham}_0(M,L)),
$
and let \(p\in L\). Choose capping maps
$
u_{\gamma(\{\psi^{(0)}_{t}\},p)},
\,
u_{\gamma(\{\psi^{(1)}_{t}\},p)}
\colon (D^2,\partial D^2)\longrightarrow(M,L)
$
whose boundaries are mapped to
$
C(\{\psi^{(0)}_{t}\},p)$ and $
C(\{\psi^{(1)}_{t}\},p),
$
respectively.

By standard homotopy theory
$
u_{\gamma(\{\psi^{(0)}_{t}\},p)}
\#
u_{\gamma(\{\psi^{(1)}_{t}\},p)}
$
is a capping of
$
C(\{\psi^{(0)}_{t}*\psi^{(1)}_{t}\},p).
$
Moreover,
\[
\begin{aligned}
\int_{D^{2}}
\left(
u_{\gamma(\{\psi^{(1)}_{t}\},p)}
\#
u_{\gamma(\{\psi^{(2)}_{t}\},p)}
\right)^*(\omega)
&=
\int_{D^{2}}
u_{\gamma(\{\psi^{(1)}_{t}\},p)}^*(\omega)
\\
&\quad+
\int_{D^{2}}
u_{\gamma(\{\psi^{(2)}_{t}\},p)}^*(\omega).
\end{aligned}
\]
and 
\[
\begin{aligned}
\mathcal{A}^{(L,\alpha)}
\bigl([\psi^{(0)}_{t}*\psi^{(1)}_{t}]\bigr)
&=
\mathcal{A}^{(L,\alpha)}
\bigl([\psi^{(0)}_{t}]\bigr)
+
\mathcal{A}^{(L,\alpha)}
\bigl([\psi^{(1)}_{t}]\bigr),
\end{aligned}
\]
in
\(\mathbb{R}/\mathcal{P}_{2}(M,\omega)\).
Hence,
\(\mathcal{A}^{(L,\alpha)}\) is a group homomorphism.
\end{proof}

We will henceforth refer to
\[
\mathcal{A}^{(L,\alpha)}\colon
\pi_{1}(\operatorname{Ham}_0(M,L))
\longrightarrow
\mathbb{R}/\mathcal{P}_{2}(M,\omega)
\]
as \emph{Weinstein's relative morphism} associated with the oriented
Lagrangian submanifold \(L\) equipped with the volume form \(\alpha\).

\medskip

The main motivation for defining the morphisms
\(\mathcal{A}_{1}^{(L,\alpha)}\)
is to detect when two Lagrangian submanifolds are not Hamiltonian isotopic.
In the following proposition, we exhibit a situation in which Weinstein's
relative morphism captures information that is not detected solely by the
global topology of the group \(\operatorname{Ham}_0(M,L)\).

More precisely, we establish a naturality identity. Conjugation by any
\(\phi\in\operatorname{Ham}(M,\omega)\) induces a group isomorphism
\[
\Phi_\phi\colon
\pi_{1}(\operatorname{Ham}_0(M,L))
\longrightarrow
\pi_{1}(\operatorname{Ham}_0(M,\phi(L))),
\]
defined by
$
\Phi_\phi\bigl([\psi_{t}]\bigr)
:=
[\phi\circ\psi_{t}\circ\phi^{-1}].
$
The following result establishes the naturality of the relative Weinstein
morphism:
\[
\mathcal{A}^{(L,\alpha)}
=
\mathcal{A}^{(\phi(L),\alpha)}
\circ\Phi_\phi.
\]

\begin{prop}
\label{p:samevalue}
Let \(L\subset(M,\omega)\) be a Lagrangian submanifold, and let
\[
[\psi_{t}]\in
\pi_{1}(\operatorname{Ham}_0(M,L)).
\]
Then, for any \(\phi\in\operatorname{Ham}(M,\omega)\) and any volume form
\(\alpha\) on \(L\),
\[
\mathcal{A}^{(L,\alpha)}
\bigl([\psi_{t}]\bigr)
=
\mathcal{A}^{(\phi(L),\alpha)}
\bigl([\phi\circ\psi_{t}\circ\phi^{-1}]\bigr)
\]
in $\mathbb{R}/\mathcal{P}_{2}(M,\omega).$
\end{prop}

\begin{proof}
Fix \( \phi \in \textup{Ham}(M,\omega) \), and let \( p \in L \). Suppose 
$
u_{ \gamma(\{\psi_{t}\},p) }
\colon (D^2,\partial D^2)\longrightarrow(M,L)
$
is a capping map for 
$
C(\{\psi_{t}\}, p) \subset L.
$
Then the composition
$
\phi \circ 
u_{ \gamma(\{\psi_{t}\},p) } \colon (D^2,\partial D^2)\longrightarrow(M,L)
$
defines a capping map for
$
C\big(\{\phi \circ \psi_{t} \circ \phi^{-1}\}, \phi(p)\big) \subset \phi(L),
$
and both maps have the same symplectic area.
Conversely, if 
$
v \colon D^{2} \to M
$
is a capping map for 
$
C\big(\{\phi \circ \psi_{t} \circ \phi^{-1}\}, q\big) \subset \phi(L),
$
then \( \phi^{-1} \circ v \) gives a capping map for 
$
C(\{\psi_{t}\}, \phi^{-1}(q)) \subset L,
$
again with equal symplectic area.

Therefore, averaging over \( L \) and \( \phi(L) \) using the same volume form \( \alpha \), we conclude that
\[
\mathcal{A}^{(L, \alpha)}([\psi_{t}]) 
= \mathcal{A}^{(\phi(L), \alpha)}([\phi \circ \psi_{t} \circ \phi^{-1}]). \qedhere
\]
\end{proof}

Thus if
  $\phi\in\operatorname{Ham}(M,\omega)$
satisfies $\phi(L_0)=L_1$, then  Proposition~\ref{p:samevalue} states that
 the two morphisms, $\mathcal{A}^{(L_0,\alpha)}$ and $\mathcal{A}^{(L_1,\alpha)}$,
 have the same image.
Hence the image of the relative Weinstein morphism, as a subgroup of
$\mathbb{R}/\mathcal{P}_2(M,\omega)$, is an invariant of the Hamiltonian
isotopy class of $L$, and it gives a criterion for deciding when two
Lagrangian submanifolds are not Hamiltonian isotopic.

\begin{cor}
\label{c:image-invariant}
Let $(M,\omega)$ be a closed, simply connected symplectic manifold and
$L_0,L_1\subset(M,\omega)$ be oriented Lagrangian submanifolds
which are difeomorphic
 with a volume form $\alpha$. If
\[
\operatorname{Im}\mathcal{A}^{(L_0,\alpha)}
\neq
\operatorname{Im}\mathcal{A}^{(L_1,\alpha)}
\quad\text{in } \mathbb{R}/\mathcal{P}_{2}(M,\omega),
\]
as subgroups of $\mathbb{R}/\mathcal{P}_2(M,\omega)$, then $L_0$ and $L_1$
are not Hamiltonian isotopic.
\end{cor}

We apply this result to the Lagrangian toric fibers of symplectic toric manifolds;
 see Sections ~\ref{sec:complexproj} and \ref{sec:toric}.

There is a special case in which Weinstein's relative morphism admits a
particularly simple and explicit description. Namely, suppose that the
traces of all points of \(L\) under the family
\(\{\psi_{t}\}\) are homotopic in \(L\). In this case, the
\(\omega\)-area of a capping of the trace is independent of the choice of
the point in \(L\), and hence the relative Weinstein morphism is simply given
by this common value.

\begin{prop}
\label{p:homotopic}
Suppose that for every pair of points $p, q \in L$, the traces 
$C(\{\psi_{t}\}, p)$ and $C(\{\psi_{t}\}, q)$ are homotopic in $L$.
Then the value
\[
\int_{D^{2}} u_{ \gamma(\{\psi_{t}\},p) } ^*(\omega) \in \mathbb{R}/ \mathcal{P}_{2}(M, \omega)
\]
is independent of the point $p \in L$.
Moreover, if $A \in \mathbb{R}/ \mathcal{P}_{2}(M, \omega)$ denotes this common value, then
\[
\mathcal{A}_{1}^{(L, \alpha)} ( [\psi_{t}]) = A.
\]
\end{prop}

\begin{proof}
Let
\[
\boldsymbol{\Psi}\colon [0,1]\times S^{1}\longrightarrow L
\]
be a homotopy between $C(\{\psi_t\},p)$ and $C(\{\psi_t\},q).$
Consider the corresponding caps
$
u_{\gamma(\{\psi_t\},p)},
u_{\gamma(\{\psi_t\},q)}
\colon (D^2,\partial D^2)\longrightarrow(M,L).
$
The images      
$
u_{\gamma(\{\psi_t\},p)}$ and
$u_{\gamma(\{\psi_t\},q)}$
 together with
$
\boldsymbol{\Psi}\bigl([0,1]\times S^{1}\bigr),
$
up to the appropriate orientation,
form a $2$-sphere in $M$. 
Since
$
\boldsymbol{\Psi}\bigl([0,1]\times S^{1}\bigr)\subset L
$
we have that
\[
\int_{D^{2}}
u_{\gamma(\{\psi_t\},p)}^*(\omega)
-
\int_{D^{2}}
u_{\gamma(\{\psi_t\},q)}^*(\omega)
\]
belongs to $\mathcal{P}_{2}(M,\omega)$.
Thus,  the value of the \(\omega\)-area of
a capping of \(C(\{\psi_{t}\},p)\)
in
\(\mathbb{R}/\mathcal{P}_{2}(M,\omega)\)
 is independent of the choice of
\(p\in L\).

Since
$
\mathcal{A}_{1}^{(L,\alpha)}([\psi_t])
$
is defined as the average of these values over $L$, the result is independent of the choice 
of $p\in L$. This proves the second assertion of the proposition.
\end{proof}

We now give an alternative description of $\mathcal{A}^{(L,\alpha)}$ when
$L$ is connected. This description facilitates the computation of the
morphism.

\begin{prop}
\label{p:factor}
Let $L\subset(M,\omega)$ be a connected Lagrangian submanifold, let
$[\psi_t]\in\pi_1(\operatorname{Ham}_0(M,L))$, and let $p\in L$. If
$u\colon(D^2,\partial D^2)\to(M,L)$ is a capping of
$\gamma(\{\psi_t\},p)$, then
\[
\mathcal{A}^{(L,\alpha)}([\psi_t])=\int_{D^2}u^*\omega
\quad\text{in }\mathbb{R}/\mathcal{P}_2(M,\omega).
\]
In particular,
$\operatorname{Im}\mathcal{A}^{(L,\alpha)}\subseteq
\omega\bigl(\pi_2(M,L)\bigr)/\mathcal{P}_2(M,\omega)$.
\end{prop}
\begin{proof}
The proof is similar to the previous result.
Let $\gamma\colon[0,1]\to L$ be a path joining $p$ to $q$. Since each
$\psi_t$ preserves $L$, the map
\[
F\colon[0,1]\times[0,1]\longrightarrow L,\qquad F(s,t)=\psi_t(\gamma(s)),
\]
defines a cylinder in the Lagrangian $L$.
If $
u_{\gamma(\{\psi_t\},p)},\,u_{\gamma(\{\psi_t\},q)}\colon
(D^2,\partial D^2)\longrightarrow(M,L)
$
are arbitrary  capping maps of the loops $\{\psi_t(p)\}$ and $\{\psi_t(q)\}$.
Then the images of $u_{\gamma(\{\psi_t\},p)},\,u_{\gamma(\{\psi_t\},q)}$ and $F$
define a sphere in $M$ and 
\[
\int_{D^{2}}
u_{\gamma(\{\psi_t\},p)}^*(\omega)
=
\int_{D^{2}}
u_{\gamma(\{\psi_t\},q)}^*(\omega)
\]
in $\mathbb{R}/\mathcal{P}_2(M,\omega).$
\end{proof}

As a consequence of Proposition~\ref{p:factor}, the definition of
$\mathcal{A}$ extends to the case where $L\subset(M,\omega)$ is a connected,
not necessarily orientable, Lagrangian submanifold. In this case, for
$[\psi_t]\in\pi_1(\operatorname{Ham}_0(M,L))$, we set
\[
\mathcal{A}^{L}([\psi_t]):=\int_{D^2}u^*\omega
\quad\text{in }\mathbb{R}/\mathcal{P}_2(M,\omega),
\]
where $u\colon(D^2,\partial D^2)\to(M,L)$ is any capping of the trace
$\gamma(\{\psi_t\},p)$ of a point $p\in L$. By
Proposition~\ref{p:factor}, this value does not depend on the choice of
$p\in L$ or of the capping, and when $L$ is oriented it agrees with
$\mathcal{A}^{(L,\alpha)}([\psi_t])$ for every volume form $\alpha$ on $L$.

\section{ A description of  $\mathcal{A}^{(L,\alpha)}$ using the Hamiltonian function} 

There is a special case in which Weinstein's relative morphism admits a
particularly simple and explicit description. In this case, the morphism can
be expressed directly in terms of the Hamiltonian generating the corresponding
loop, providing a useful alternative description of the invariant and
clarifying its relation to the original Weinstein morphism.

In this section, we focus on the group
\(\pi_1(\operatorname{Ham}_0(M,L))\).
Let \(\{\psi_t\}_{t\in[0,1]}\) be a loop of Hamiltonian diffeomorphisms
generated by the time-dependent Hamiltonian
\(H_t\colon M\to\mathbb{R}\). We assume that \(H_t\) is normalized for every
\(t\), that is,
\[
\int_M H_t\,\omega^n=0.
\]
Recall that an alternative expression for Weinstein's morphism
\[
\mathcal{A}^M\colon
\pi_1(\operatorname{Ham}(M,\omega))
\longrightarrow
\mathbb{R}/\mathcal{P}_2(M,\omega)
\]
is given by
\[
\mathcal{A}^M([\psi_t])
=
\int_{D^2}
u_{\gamma(\{\psi_t\},p_0)}^*(\omega)
-
\int_0^1 H_t(\psi_t(p_0))\,dt,
\]
where
\(u_{\gamma(\{\psi_t\},p_0)} \colon (D^2,\partial D^2)\longrightarrow(M,L) \)
is a filling map as above and \(p_0\in M\) is a fixed base point.

We next derive an analogous expression for the relative Weinstein morphism.
The resulting formula relates the relative morphism to the original Weinstein
morphism and expresses the difference between the two in terms of the
Hamiltonian generating the loop and the Lagrangian submanifold \(L\).

\begin{thm}
\label{t:wholevsrel}
Let \(\{\psi_t\}_{t\in[0,1]}\) be a loop in
\(\operatorname{Ham}_0(M,L)\) based at the identity and generated by the
normalized Hamiltonian \(H_t\). If \(L\) is connected and
\(p_0\in L\) is a fixed base point, then
\begin{eqnarray*}
\mathcal{A}_{1}^{(L,\alpha)}([\psi_t])
=
\mathcal{A}^M([\psi_t])
+
\int_0^1 H_t(\psi_t(p_0))\,dt 
\in \mathbb{R}/\mathcal{P}_2(M,\omega).
\end{eqnarray*}
\end{thm}

\begin{proof}
As discussed in the definition of
\(\mathcal{A}_{1}^{(L,\alpha)}\), let
$
F_{\psi_t}\colon M\longrightarrow\mathbb{R}
$
be a lift of
$
A\circ\operatorname{ev}_{\psi_t}\colon
M\longrightarrow
\mathbb{R}/\mathcal{P}_2(M,\omega).
$
Before passing to the quotient by the period group, and with a slight abuse
of notation, we may regard the following equality as an equality in
\(\mathbb{R}\),
 \[
\mathcal{A}^M([\psi_t])
=
F_{\psi_t}(p_0)
-
\int_0^1 H_t(\psi_t(p_0))\,dt \in \mathbb{R}.
\]
Since \(\{\psi_t\}_{t\in[0,1]}\) is a loop in
\(\operatorname{Ham}_0(M,L)\), each \(\psi_t\) preserves \(L\), and hence
the Hamiltonian vector field \(X^H_{t}\) is tangent to \(L\). Since \(L\subset (M,\omega)\) is
Lagrangian and connected, then \(H_t\) is constant on \(L\) for every \(t\). 

Upon integrating the preceding expression over \(L\), we obtain
\begin{eqnarray*}
{\operatorname{Vol}(L,\alpha)} \cdot \mathcal{A}^M([\psi_t])
&=&
\int_ L 
F_{\psi_t}(p_0) \, \alpha
-
{\operatorname{Vol}(L,\alpha)} \cdot\int_0^1 H_t(\psi_t(p_0))\,dt 
\end{eqnarray*}
and the proposition follows.
\end{proof}

\section{Example: The 2-sphere}
\label{s:sphere}

Consider the \(2\)-sphere \((S^2,\omega)\) of radius \(1\). Its spherical
period group is
$
\mathcal{P}_2(S^2,\omega)
=
\mathbb{Z}\langle 4\pi\rangle.
$
Let \((L,\alpha)\) be the circle of length \(2\pi\), and, for
\(h\in[0,1)\), let \(L_h\subset S^2\) be the circle in the plane \(z=h\).

The rotation about the \(z\)-axis, denoted by \(\{\psi_t\}\), preserves
\(L_h\subset  (S^2,\omega)  \). Thus, for every \(h\in[0,1)\), it determines an element
$
[\psi_t]\in\pi_1(\operatorname{Ham}(S^2,L_h)).
$
The \(\omega\)-area of the northern cap bounded by \(L_h\) is
$
2\pi(1-h),
$
while the trace of any point of \(L_h\) under \(\{\psi_t\}\) is precisely
\(L_h\) itself. Therefore, by Proposition~\ref{p:homotopic},
\[
\mathcal{A}^{(L_h,\alpha)}[\psi_t]
=
2\pi(1-h)
\in
\mathbb{R}/\mathbb{Z}\langle4\pi\rangle.
\]

Since $\pi_2(S^2, L_h )$ is isomotphic to $\pi_2(S^2)\oplus \pi_1(L_h)$, 
it follows that $\textup{Im} (\mathcal{A}^{(L_h,\alpha)})$ is $\mathbb{Z}\langle 2\pi(1-h) \rangle /  
\mathbb{Z}\langle4\pi\rangle.$ Thus, by
Corollary \ref{c:image-invariant}
 for distinct
\(h_0,h_1\in[0,1)\), the Lagrangian submanifolds
\(L_{h_0}\) and \(L_{h_1}\) are not Hamiltonian isotopic, provided their associated cyclic groups in $\mathbb{R}/\mathbb{Z}$ are distinct. This recovers, in
this particular case, the well-known fact that the circles \(L_h\) are
pairwise non-Hamiltonian-isotopic.

Moreover, if \(h\in[0,1)\) is such that
$(1-h)/{2}$ is irrational,
then
$
\mathcal{A}^{(L_{h},\alpha)}_1[\psi_t]
=
2\pi(1-h)
$
has infinite order in
\(\mathbb{R}/\mathbb{Z}\langle4\pi\rangle\). Since
\(\mathcal{A}^{(L_{h},\alpha)}_1\) is a group homomorphism, it follows that
\[
[\psi_t]\in
\pi_1(\operatorname{Ham}(S^2,L_{h}))
\]
has infinite order.
This completes the proof of Theorem~\ref{t:sphere} stated in the
Introduction.

 \section{Example: Lagrangian torus fibers in \(\mathbb{C}P^n\)}
 \label{sec:complexproj}

 We next consider a nontrivial example arising from the toric structure of
 complex projective space. This example is particularly useful because the
 relative Weinstein morphism can be computed explicitly on the standard
 loops generated by the toric action, and the resulting values depend
 directly on the moment-map coordinates of the Lagrangian torus.

 Consider \(\mathbb{C}P^n\) equipped with its standard toric symplectic
 structure. We normalize the Fubini--Study symplectic form so that the
 symplectic area of a complex projective line is \(1\). The moment polytope
 is then the standard simplex
 $
 \Delta (\mathbb{C}P^n,\omega_{\mathrm{FS}})
 =
 \operatorname{conv}
 \left\{
 (0,\ldots,0),
 (1,0,\ldots,0),
 \ldots,
 (0,\ldots,0,1)
 \right\}
 \subset \mathfrak{t}^*\simeq \mathbb{R}^n.
 $
 and 
 $
 \mathcal{P}_2(\mathbb{C}P^n,\omega_{\mathrm{FS}})
 =
 \mathbb{Z}.
 $
 Let
 $
 \mathbf a=(a_1,\ldots,a_n)\in\operatorname{int}(\Delta (\mathbb{C}P^n,\omega_{\mathrm{FS}}) ),
 $ be an interior point
 and let
 $
 L_{\mathbf a}:=\mu^{-1}(\mathbf a)
 $
 be the corresponding Lagrangian torus. Fix a volume form
 \(\alpha\) on \(L_{\mathbf a}\). We will compute
 \[
 \mathcal{A}^{(L_{\mathbf a},\alpha)}
 \colon
 \pi_1\bigl(\operatorname{Ham}_0(\mathbb{C}P^n,L_{\mathbf a})\bigr)
 \longrightarrow
 \mathbb{R}/\mathbb{Z}
 \]
 on the standard circle subgroups of the toric action.

 The torus \(\mathbb{T}^n\) acts on \(\mathbb{C}P^n\) by
 \[
 (e^{2\pi i t_1},\ldots,e^{2\pi i t_n})
 \cdot
 [z_0,\ldots,z_n]
 =
 [z_0,e^{2\pi i t_1}z_1,\ldots,e^{2\pi i t_n}z_n],
 \]
 and the corresponding  moment map
 $
 \mu\colon \mathbb{C}P^n\longrightarrow\mathbb{R}^n
 $
 is given by
 $$
 \mu([z_0,\ldots,z_n])
 =
 \left(
 \frac{|z_1|^2}{\|z\|^2},
 \ldots,
 \frac{|z_n|^2}{\|z\|^2}
 \right).
 $$
 For \(1\leq k\leq n\), consider the \(k\)-th circle subgroup
 \[
 \psi^{(k)}_t
 [z_0,\ldots,z_n]
 =
 [z_0,\ldots,e^{2\pi i t}z_k,\ldots,z_n].
 \]
 Then  its normalized Hamiltonian is
 \[
 H^{(k)}([z_0,\ldots,z_n])
 =
 \frac{|z_k|^2}{\|z\|^2}
 -\frac{1}{n+1}.
 \]
 The toric action preserves every moment-map fiber, therefore
 $
 \{\psi^{(k)}_t\}$ induces an element in 
 $\pi_1(
 \operatorname{Ham}_0(\mathbb{C}P^n,L_{\mathbf a})).
 $
 Note that for every \(p\in L_{\mathbf a}\), we have
 \[
 H^{(k)}(\psi^{(k)}_t(p))
 =
 a_k-\frac{1}{n+1}.
 \]
 On the other hand, the absolute Weinstein morphism of this circle action
 is
 \[
 \mathcal{A}^{\mathbb{C}P^n}_1
 ([\psi^{(k)}_t])
 =
 \frac{1}{n+1} \in \mathbb{R}/\mathbb{Z}.
 \]
 Therefore, by Theorem~\ref{t:wholevsrel},
 \[
 \begin{aligned}
 \mathcal{A}^{(L_{\mathbf a},\alpha)}_1
 ([\psi^{(k)}_t])
 &=
 \mathcal{A}^{\mathbb{C}P^n}_1
 [\psi^{(k)}_t]
 +
 H^{(k)}(\psi^{(k)}_t(p_0))
 \\
 &=
 \frac{1}{n+1}
 +
 a_k-\frac{1}{n+1}
 \\
 &=
 a_k  \in \mathbb{R}/\mathbb{Z}.
 \end{aligned}
 \]

 These calculations prove Theorem~\ref{t:aoncpn}. 
 Now $\pi_2(\mathbb{C}P^n, L_{\mathbf a})$ is isomorphic to
 $\pi_2(\mathbb{C}P^n)\oplus \pi_1( L_{\mathbf a})$, thus
 $\textup{Im} (\mathcal{A}^{(L_{\mathbf a},\alpha)})$ is $\mathbb{Z}\langle a_1,\ldots,a_n \rangle /  
 \mathbb{Z}.$
 In particular, 
 for $
 \mathbf a=(a_1,\ldots,a_n), \mathbf b=(b_1,\ldots,b_n)
 \in\operatorname{int}(\Delta (\mathbb{C}P^n,\omega_{\mathrm{FS}}) ),
 $
 the Lagrangian tori
 \[
 L_{\mathbf a}=\mu^{-1}(\mathbf a),
 \qquad
 L_{\mathbf b}=\mu^{-1}(\mathbf b)
 \]
 are not Hamiltonian isotopic whenever the groups
  $\mathbb{Z}\langle a_1,\ldots,a_n \rangle$ and  $\mathbb{Z}\langle b_1,\ldots,b_n \rangle$
 are distinct.

\section{Example: Lagrangian torus fibers in toric manifolds}
\label{sec:toric}

In this section we calculate the value of the morphism on the standard rotations generated by the toric action on a toric manifold. More precisely, let $(M^{2n},\omega)$ be a compact connected symplectic manifold equipped with an effective Hamiltonian action of the torus $\mathbb{T}^n$. Let $\mathfrak{t}$ denote the Lie algebra of $\mathbb{T}^n$; we have the identification $\mathfrak{t}\cong \mathbb{R}^n$. For $\xi \in \mathfrak{t}$, the moment map $\mu : M \rightarrow \mathfrak{t}^{*}$ defines the Hamiltonian function $H_{\xi}:M \rightarrow \mathbb{R}$ defined by \[p\mapsto \langle \mu,\xi\rangle(p) = \mu(p)(\xi).\]The Hamiltonian vector field generated by this function is given by \[X_{\xi}(p) = \dfrac{d}{dt}\Big|_{t=0}\exp(t\xi)\cdot p.\] 

The image of $M$ under the moment map, denoted by $\Delta = \mu(M)$, is a Delzant polytope. We gather some facts about the preimages of points in the $\Delta(M,\omega)$:
\begin{itemize}
\item For ${\bf a}\in\text{Int} (\Delta(M,\omega))$, the Lagrangian torus $L_{\bf a} = \mu^{-1}({\bf a})$ corresponds to a torus $T^n$. 
    \item For a vertex $v$ in the polytope $\Delta$, the preimage correspond to a single point $\mu^{-1}(v) = \{x_v\}$ and this point is preserved by the Hamiltonian action of $\mathbb{T}^n$.
\end{itemize}
Under the identifications $\mathfrak{t}\equiv\mathbb{R}^n$ and $\mathbb{T}^n \equiv \mathbb{R}^n / \mathbb{Z}^n$, the integral lattice can be identified by $\Lambda \equiv \mathbb{Z}^n$. For $\xi \in \Lambda$ and $t\in[0,1]$ define the path \[\psi^{\xi}_t:M\rightarrow M\] by $\psi^{\xi}(p) = \exp(t\xi)\cdot p$. Note that since $\exp(\xi) = Id$, this path corresponds to a loop at the of Hamiltonian diffeomorphisms that preserve the Lagrangians $L_{\bf a}$, \[\psi^{\xi} \in \pi_1\left(\operatorname{Ham}(M,L_{\bf a})\right). \]

The strategy to compute $\mathcal{A}_{1}^{(L_{\bf a},\alpha)}([\psi^{\xi}])$ is to find and explicit caping disk of the loop $\psi^{\xi}(p)$ for $p\in L_{\bf a}$. For a vertex $v \in \Delta(M,\omega)$, let $x_v$ the preimage of the vertex under the moment map. Since $M$ is connected we can find a path $c:[0,1] \rightarrow M$ such that
\begin{eqnarray*}
    c(0) = x_v, && c(1) = p,
\end{eqnarray*}
and we require that $c$ remains constant near zero. We apply $\psi^{\xi}$ to the path $c(t)$ to obtain a disk $B(s,t) = \psi_{t}^{\xi}(c(s))$ (see Fig.\ref{fig:capping_disk} below) 

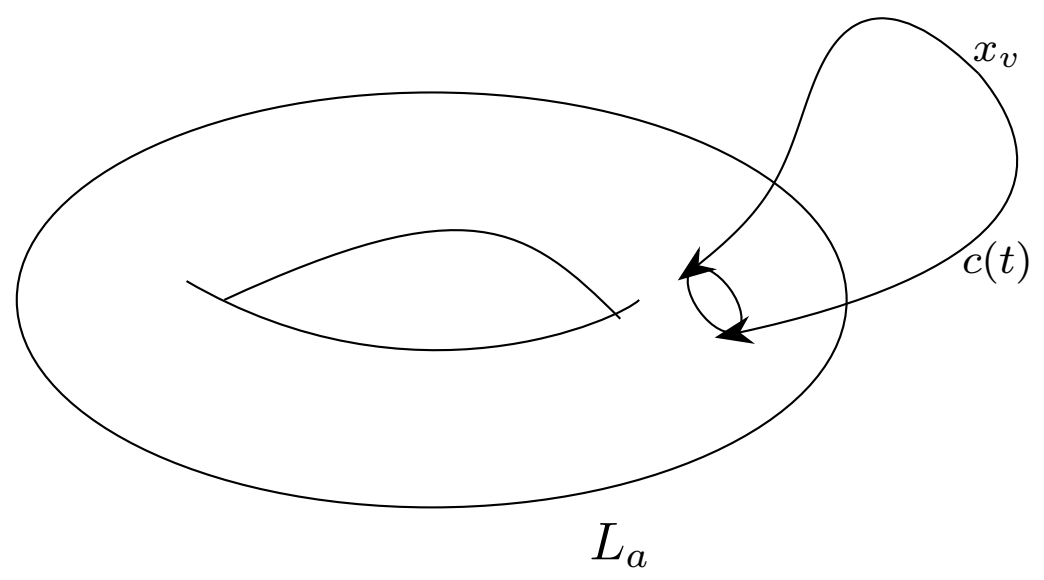
\begin{figure}[!ht]
    \centering
    \resizebox{0.5\textwidth}{!}{%
    \begin{circuitikz}

    \tikzstyle{every node}=[font=\fontsize{10.2pt}{13.3pt}\selectfont]

    \draw (5.75,11.5) ellipse (2.75cm and 1.375cm);

    \draw (4.125,11.625) .. controls (5.625,10.75) 
        and (7,11.375) .. (7.125,11.5);

    \draw (4.375,11.5) .. controls (6,12.25) 
        and (6.375,12) .. (7,11.375);

    \draw [-{Stealth[scale=1.5]}]
        (9.375,13) .. controls (10,12.25) 
        and (9.5,11.625) .. (7.625,11.25);

    \node [
        font=\fontsize{8.0pt}{10.4pt}\selectfont,
        inner xsep=0.080cm,
        inner ysep=0.085cm
    ] at (9.5,13.125) {$x_v$};

    \node [
        font=\fontsize{8.0pt}{10.4pt}\selectfont,
        inner xsep=0.080cm,
        inner ysep=0.085cm
    ] at (9.5,11.75) {$c(t)$};

    \draw [rotate around={-54:(7.625,11.5)}]
        (7.625,11.5) ellipse (0.25cm and 0.125cm);

    \draw [-{Stealth[scale=1.5]}]
        plot[variable=\t, domain=0:1, samples=100]
        (
        {(1*(1-\t)^4)*(9.375)
        + (4*(1-\t)^3*\t^1)*(8.25)
        + (6*(1-\t)^2*\t^2)*(8.125)
        + (4*(1-\t)^1*\t^3)*(8.375)
        + (1*\t^4)*(7.375)},
        {(1*(1-\t)^4)*(13)
        + (4*(1-\t)^3*\t^1)*(14.125)
        + (6*(1-\t)^2*\t^2)*(12.5)
        + (4*(1-\t)^1*\t^3)*(12.375)
        + (1*\t^4)*(11.625)}
        );

    \node [
        font=\fontsize{10.2pt}{13.3pt}\selectfont,
        inner xsep=0.080cm,
        inner ysep=0.085cm
    ] at (7,9.875) {$L_a$};

    \end{circuitikz}%
    }
    \caption{Capping disk for $\psi^{\xi}$}
    \label{fig:capping_disk}
\end{figure}


\noindent Recall the Hamiltonian function $H^{\xi}(c(s)) = \langle \mu(c(s)),\xi\rangle$ generates the vector field $X_\xi (B(s,t)) = \partial_t B$. Therefore, 

\begin{eqnarray*}
    B^{*}\omega(\partial_s ,\partial_t) &=& \omega(\partial_s B, \partial_t B) \\
    &=& \omega(\partial_s B, X_{\xi}(B(s,tp))) \\
    &=& \dfrac{d}{ds}H^{\xi}(c(s)) \\
    &=& \dfrac{d}{ds} \langle \mu(c(s)), \xi \rangle.
\end{eqnarray*}

\noindent The area of the capping disk $B(s,t)$ is then computed as follows

\begin{eqnarray*}
    \int_{D^2}B^{*}\omega &=& \int_{0}^{1} \int_{0}^{1} \dfrac{d}{ds} \langle \mu(c(s)), \xi \rangle \\
                        &=& \langle \mu(p), \xi \rangle - \langle \mu (x_v), \xi \rangle \\
                        &=& \langle {\bf a} - v, \xi \rangle,
\end{eqnarray*}
where we used that $p\in L_a$. That is

\begin{equation}\label{eqn:toric_formula}
    \mathcal{A}_{1}^{(L_{\bf a},\omega)}([\psi^{\xi}]) = \langle a - v, \xi \rangle
\end{equation}

\noindent In general, we can describe a Delzant polytope as \[\Delta (M, \omega)= \{ x\in\mathfrak{t}^* : f_i(x) \geq 0, i = 1,\cdots,d\}\] and the facets are determined by $F_i = \{f_i = 0\}$ for $i=1,\cdots,d$. In general the defining equation of a facet can be written as follows: 
\[f_i(x) = \langle x, \eta_i\rangle - \lambda_i,\]
where $\eta_i$ is the normal primitive inward i.e the minimal generator of the rotation in the $\eta_i$ direction and $\lambda_i$ the affine constant. Consider a path $c(s)$ connecting a point $p\in L_a$ and a point $x_{v_i}$ that is the preimage of a vertex $v_i$ in the $i$-th facet. The following results shows that the morphism $\mathcal{A}^{(L_{\bf a},\omega)}$ on the primitive rotations can be understood in terms of the affine distance between an interior point to the facet. 

\begin{prop}
    For a primitive inward rotation $\eta_i \in \Lambda$, we have  \[\mathcal{A}^{(L_{\bf a},\omega)}([\psi^{\eta_i}])= f_i({\bf a}). \]   
\end{prop} 
\begin{proof}
    Equation (\ref{eqn:toric_formula}) and the fact that $v_i$ is in the $i$-th facet imply 
    \begin{eqnarray*}
    \mathcal{A}^{(L_{\bf a},\omega)}([\psi^{\eta_i}]) &=& \langle {\bf a}-v_i, \eta_i\rangle\\
    &=& \langle {\bf a}, \eta_i \rangle - \langle v_i , \eta_i\rangle \\
    &=& \langle {\bf a}, \eta_i \rangle - \lambda_i \\
    &=& f_i({\bf a}). 
\end{eqnarray*}
\end{proof}

\noindent The toric action defines homomorphisms
\[\iota_{\bf a} : \mathbb{Z}^n \rightarrow \pi_1 \left(\textup{Ham}_0(M,L_{\bf a})\right),\]
and 
\[\iota : \mathbb{Z}^n \rightarrow \pi_1 \left(\textup{Ham}(M,\omega)\right).\] Consider the inclusion \[j: \textup{Ham}_0(M,L_{\bf a} ) \hookrightarrow \textup{Ham}(M,\omega),\]

\begin{thm} For ${\bf a}\in \textup{Int}(\Delta)(M,\omega)$ we have
    \begin{enumerate}
        \item The homomorphism $\iota_{\bf a} : \mathbb{Z}^n \rightarrow \pi_1 \left(\textup{Ham}_0(M,L_{\bf a})\right)$ is injective. In particular, every non-zero integral toric rotation defines an element of infinite order in relative fundamental group.
        \item Let $K_\iota = \ker(\iota) \subset \mathbb{Z}^n$, then \[\iota_{\bf a} (K_\iota) = \textup{im}(\iota_{\bf a})\cap \ker(j_*). \]
        \item Define \[
        \mathcal{S} = \bigcup_{\xi \in K_\iota\setminus\{0\}, A \in \Gamma_{\omega}} \{ {\bf a}\in \textup{Int}(\Delta(M,\omega)): \langle {\bf a}-v,\xi\rangle = A\}.
        \]
        Then $\mathcal{S}$ is a countable union of proper affine hyperplanes sections of $\textup{Int}(\Delta(M,\omega))$. Moreover, for every $a\in\textup{Int}(\Delta(M,\omega)  )\setminus\mathcal{S}$, the restriction \[\mathcal{A}^{(L_{\bf a},\omega)}\Big|_{K_\iota} : K_\iota \rightarrow \mathbb{R}/\Gamma_\omega,\] is injective. In particular, the are homomorphism detects every non-trivial absloutely contractible toric rotation on a generic Lagrangian tori.
    \end{enumerate}
\end{thm}
\begin{proof}

Proof of (1): Let ${\bf a}\in\text{Int}(\Delta(M,\omega)   )$ and fix $p\in L_{\bf a}$. Since ${\bf a}$ is in the interior, the torus action on the Lagrangian $L_{\bf a}$ is free and transitive. Therefore the map $f_p : \mathbb{T}^n \rightarrow L_a$, $f_p (g) = g\cdot p$ is a difeomorphism. On the other hand, the map $f_p$ factors as 
\[
\mathbb{T}^n \hookrightarrow \textup{Ham}_0(M,L_{\bf a}) \xrightarrow{ev_p} L_{\bf a}.
\]
Taking fundamental groups we obtain
\[\mathbb{Z}^n \xrightarrow{\iota_{\bf a}} \pi_1\left(\textup{Ham}_0(M,L_{\bf a}) \right)\xrightarrow{(ev_p)_*} \pi_1\left(L_{\bf a}\right),\]
the composition agrees with $(f_p)_*$ ; after identifying \[
\pi_1(L_a)\cong\pi_1(T^n)\cong\Lambda\cong\mathbb Z^n.
\]we get \[(ev_p)_* \circ \iota_{\bf a} = \text{Id}_{\mathbb{Z}^n},\] proving (1). In fact, this statement proves more. Since the fundamental group of a topological group is abelian, we obtain th splitting 
\[\pi_1 \left( \textup{Ham}_0(M,L_{\bf a}) \right) \cong \mathbb{Z}^n \oplus \ker((ev_p)_*).\] 

Proof of (2): Let $K_\iota = \ker(\iota)$, by definition  \[\iota = j_* \circ \iota_{\bf a}.\] Hence, $\xi\in K_\iota$, if and only if $\iota_{\bf a}(\xi) \in \ker(j_*).$ It follows that \[\iota_{\bf a} (K_\iota) = \text{im}(\iota_{\bf a}) \cap \ker(j_*).\] Note that $K_\iota$ can be regarded as subgroup of $\mathbb{Z}^n$ and so it itself is free abelian, say $K_\iota \cong \mathbb{Z}^r$ for $0\leq r\leq n$. In particular we have the inclusion,  \[\mathbb{Z}^r \hookrightarrow \ker(j_*).\] The rank $r$ is independent of ${\bf a}$, because $\iota$ is defined entirely by the ambient torus action in $M$.\\

\noindent Proof of (3): If $K_\iota = \{0\}$, then the result follows. Thus, we assume $K_\iota \neq \{0\}$. For each $\xi \in K_\iota \setminus \{0\}$ define \[H_\xi = \{{\bf a}\in\text{Int}(\Delta(M,\omega)) : \langle {\bf a}-v,\xi \rangle \in \Gamma_\omega. \}\] Consequently, we have the alternative description \[H_\xi = \bigcup_{A \in \Gamma_\omega} \{ {\bf a}\in \text{Int}(\Delta(M,\omega)): \langle {\bf a}-v,\xi \rangle = A\}.\] Since $H_2(M;\mathbb{Z})$ is finitely generated, $\Gamma_\omega$ is a countable subgroup of $\mathbb{R}.$ Therefore, since $K_\iota$ is also countable we have that $\mathcal{S}$ is a countable union of countable sets.

\noindent For ${\bf a}\in \text{Int}(\Delta(M,\omega)) \setminus \mathcal{S}$, for every non-zero $\xi \in K_\iota$ we have 
$\mathcal{A}^{(L_{\bf a},\alpha)}(\xi) \neq 0$. Therefore $\ker(\mathcal{A}^{(L_{\bf a},\alpha)}|_{K_\iota}) = \{0\}$ i.e it is injective.
\end{proof}

\noindent \textbf{Remark:} The last property of this theorem proves that generically the morphism detects non-trivial loops in the relative group preserving the lagrangian fibers that are contractible in the absoloute groups. \\

\begin{cor}
    Let $(M,\omega)$ be a compact symplectic toric manifold with moment polytope  
    \[
        \Delta(M,\omega)=\{x\in\mathfrak t^*:f_i(x)=\langle x,\eta_i\rangle-\lambda_i\geq0,\quad i=1,\ldots,d\}.
    \]
    If for $\mathbf{a},\mathbf{b} \in \text{Int}(\Delta(M,\omega) )$ we have that \[\mathbb{Z}\langle f_1 (\mathbf{a}), \cdots, f_d (\mathbf{a})\rangle \neq \mathbb{Z}\langle f_1 (\mathbf{b}), \cdots, f_d (\mathbf{b})\rangle\] then $L_{\mathbf{a}}$ and $L_{\mathbf{b}}$ are not Hamiltonian isotopic.
\end{cor}

\section{Example: The real projective space in \(\mathbb{C}P^n\)}
\label{sec:comp}

In this example we consider a Lagrangian submanifold which is nota a toric fiber.
In \((\mathbb{C}P^n,\omega_{\mathrm{FS}})\), consider the Hamiltonian loop
\[
\begin{aligned}
\psi_t([z_0,\ldots,z_n])
=
[&\cos(2\pi t)z_0-\sin(2\pi t)z_1,\\
&\sin(2\pi t)z_0+\cos(2\pi t)z_1,
z_2,\ldots,z_n].
\end{aligned}
\]
Assume that \(n\geq 2\). This loop preserves the Lagrangian submanifold
$
\mathbb{R}P^n\subset
(\mathbb{C}P^n,\omega_{\mathrm{FS}}).
$
Moreover, it fixes pointwise the projective subspace
\[
\{[x_0:\cdots:x_n]\in\mathbb{R}P^n:
x_0=x_1=0\}\cong\mathbb{R}P^{n-2}.
\]
In particular, the fixed-point set of the loop in \(\mathbb{R}P^n\) is nonempty. 
Notice that for $n$ even  \(\mathbb{R}P^n\) is non orientable, thus we used the
the detinition of $\mathcal{A}^{\mathbb{R}P^n}$ given by a capping disk.
Therefore, 
$
\mathcal{A}^{\mathbb{R}P^n}[\psi_t ]=0.
$


\bibliographystyle{acm}
\bibliography{Ref.bib}  
\end{document}